\documentclass[10pt]{article}
\usepackage{amsmath, amssymb, amsthm, pb-diagram, euscript, mathrsfs, textcomp}

\usepackage[usenames]{color}
\usepackage{colortbl}

\title{Complete Characterization of $\K$-Theory for $\Cst$-algebras Associated to Locally Finite Unoriented Graphs.}
\author{Nikolay Ivankov, Natalia Iyudu\\}

\theoremstyle{plain}
\newtheorem{prop}{Proposition}[section]

\newtheorem{lem}[prop]{Lemma}
\newtheorem{cor}[prop]{Corollary}
\newtheorem{thm}[prop]{Theorem}

\theoremstyle{definition}
\newtheorem{defn}[prop]{Definition}

\newtheorem{rem}[prop]{Remark}

\theoremstyle{remark}

\newcommand{\E}{\mathcal{E}}

\newcommand{\K}{\mathrm{K}}

\newcommand{\Nb}{\mathbb{N}}

\newcommand{\Zb}{\mathbb{Z}}

\newcommand{\Osm}{\mathcal{O}}

\newcommand{\al}{\alpha}
\newcommand{\be}{\beta}
\newcommand{\gm}{\gamma}

\newcommand{\ep}{\varepsilon}
\newcommand{\et}{\eta}

\newcommand{\la}{\lambda}

\newcommand{\Id}{\mathrm{Id}}

\renewcommand{\Im}{\mathrm{Im}}

\newcommand{\xto}{\xrightarrow}

\newcommand{\iy}{\infty}

\newcommand{\ol}[1]{\overline{#1}}

\newcommand*{\Cst}{\textup C^*}

\DeclareMathOperator{\coker}{coker}

\begin{document}

\maketitle

\begin{abstract}
In this paper we give a complete description of $\K$-theory groups for Cuntz-Krieger $\Cst$-algebras associated to general locally-finite (topologically connected) graphs via Bass-Hashimoto operator. Our result generalizes the one obtained by the second author in \cite{Iyu13} for the case of graphs with not necessarily finite first Betti numbers. On the bases of  purely graph-theoretical method introduced in \cite{CLM08} and developed further in \cite{Iyu13}, we prove that for the algebra $\Osm_\E$ associated an infinite graph $\E$ of the above form holds $\K_0(\Osm_E)=\Zb^{\beta(\E)}\oplus \Zb^{\gamma(\E)}$ and $\K_1(\Osm_\E) = \Zb^{\gamma(\E)}$, where $\beta(\E)=\dim H_1(\E)$ and $\gamma(\E)$ stand for the cardinality of the valency set of $\E$ defined in the article.
\end{abstract}

\section{Introduction}

The object of study of the present paper are locally-finite non-directed graphs, i.e.\ those that have a finite number of edges adjusted to each of its axes. 
We consider the case when these graphs have infinite number of edges (and therefore also axes), as the finite case has been studied in detail in \cite{CLM08}.


Let $\E$ be any connected locally-finite graph, possibly with loops, multiple edges and sinks. To this graph we associate a directed graph $E$ with the same set of axes and a pair of edges $\{x,\ol x\}$, pointing in both possible directions, associated to each edge of $\E$ (the loop edges get doubled). Let $E^1$ be the set of edges of $E$, and $\Zb^{|E^1|}$ be a free $\Zb$-module spanned by these edges. The \emph{Bass-Hashimoto operator} $\Phi\colon \Zb^{|E^1|} \to \Zb^{|E^1|}$ is defined on basis elements $x\in E^1$ by
\begin{equation}\label{eq:BassHashimoto}
\Phi(x)= \left(\sum_{r(x)=s(x')}x'\right) - \ol x,
\end{equation}
where the sum expands over all edges $x'\in E^1$ whose source coincides with the range of $x$. This operator has been considered by Hashimoto \cite{Has89} and Bass \cite{Bas92} in their study of Ihara zeta function of a graph.

The operator $\Phi$ associated to the graph $\E$ induces an (infinite) matrix $A_\E$ indexed by the elements of $\E^1$, which is called the incidence matrix for $\E$. To each such matrix me may apply the procedure described by Cuntz and Krieger in \cite{CK80} to define a so-called Cuntz-Krieger algebra $\Osm_\E$: this $\Cst$-algebra is generated by partial isometries $\{S_x\}_{x\in \E^1}$ subject to the relations $S_x^\ast S_x =\sum_{y\in \E^1} A_{x,y} S_y S_y^\ast$.

For finite graphs it has been shown that the $\K$-groups of an algebra $\Osm_E$ admit an algebraic description in terms of operator $\Phi$ as a map between $\Zb$-modules, namely that $K_0(\Osm_\E)\cong \coker(\Id-\Phi)$ and $\K_1(\Osm_E)\cong \ker(\Id-\Phi)$ (see e.g.\ \cite{CK80} for the case of finite incidence matrix and its generalization for infinite ones in \cite{PR96}). In \cite{CLM08} there have been established that for the algebras $\Osm_\E$ associated to finite graphs  $\E$ there exists a purely graph-theoretical description of their $\K$-groups in terms of their first Betti numbers. More precisely, it was shown that, given a finite graph $\E$ of genus not less than $2$, one has $\K_0(\Osm_\E)=\Zb^{\beta(\E)}\oplus \Zb/(\beta(\E)-1)\Zb$ and $\K_1(\Osm_\E) =\Zb^{\beta(\E)}$, where $\beta(\E)=\dim H_1(\E)$ is the first Betti number of the geometric realization of $\E$. These results were then generalized by one of the authors in \cite{Iyu13} for a particular case of infinite locally finite graphs, namely for the those with finite first Betti number. There it was proved that for such a graph $\E$ one has
\begin{equation}\label{eqn:K_0_and_K_1}
\K_0(\Osm_\E)=\Zb^{\beta(\E)}\oplus \Zb^{\gamma(\E)} \quad \text{whereas}\quad  \K_1(\Osm_\E) =\Zb^{\beta(\E)},
\end{equation}
where $\gamma(\E)$ stands for the so-called valency number of a graph. This result has also shown the inconsistency of the hypothesis that the group $\K_1(\Osm_\E)$ is the torsion-free part of $\K_0(\Osm_\E)$.

In the current paper we extend the results of \cite{Iyu13}, namely the formula \eqref{eqn:K_0_and_K_1}, to the case of general locally finite infinite graphs. Since the only case that was not studied yet is that of the graphs with ``infinite first Betti number'', some adjustments had to be made in order to render the formula \eqref{eqn:K_0_and_K_1} meaningful: thus, instead of Betti number we talk about the cardinality of the generator set of the first homology group, and introduce an invariant up to an isomorphism definition of the valency set.

The paper is organized as follows. In Section \ref{sec:Notatins} we introduce the notation that will allow us to make exposition more illustrative. Section \ref{sec:Rosetree} is devoted to the calculation of $\K$-groups for the algebras associated to a particular kind of locally finite graphs with possibly infinite first Betti numbers, which we call ``rose-tree'' graphs. Here we establish that formula \eqref{eqn:K_0_and_K_1} holds for this kind of graphs. In Section \ref{sec:TechLemmata} we generalize the Lemmata 3.2 and 4.3 of \cite{Iyu13}. These results in particular imply that there is a natural isomorphism between $\K$-groups of the algebras associated to two graphs, one of which is obtained from another by ``shrinking'' of a countable number of mutually disjoint finite tree subgraphs: an operation that preserves Betti number and valency set. In section \ref{sec:Reduction} we provide an algorithmic procedure that allows us to 
establish an isomorphism between $\K$-groups of a given graph and of a particular rose-tree graph, obtained from it. Finally, in Section \ref{sec:Cnnection} we discuss the connection of the groups $\K_0(\Osm_\E)$ and $K_1(\Osm_\E)$ to the first homology group of $\E$ and its valency set.

 Since the operator \eqref{eq:BassHashimoto} is well-defined only for locally-finite graphs,
  the result in the presented paper, combined with those of \cite{CLM08} and \cite{Iyu13} completely solve the problem of calculating the $\K$-theory for the algebras associated to these operators.





\section{Notations and Conventions}\label{sec:Notatins}


Throughout the paper we assume the initial graph $\E$ to be connected. For non-connected graphs one has $\K_i(\Osm_\E) = \bigoplus_\kappa \K_i(\Osm_{\E_\kappa})$, $i=0,1$, where $\kappa$ runs through all the connected components $\E_\kappa$ of $\E$.


\subsection{Bidirected Graphs}

By (locally finite) \emph{bidirected graph} or \emph{bi-graph} we shall understand a (locally finite) directed graph $(E^0,E^1,s,r)$ together with a bijection $\ol\cdot\colon E^1\to E^1$, such that
\begin{itemize}
\item the set of edges $E^1$ may be represented as a disjoint union
$$
E^1 = E^{1,0}\sqcup E^{1,1};
$$
\item the map $\ol\cdot$ restricts to the isomorphisms
$$
\ol\cdot\colon E^{1,i}\to E^{1,i+1},
$$
where $i=0,1$ and the sum is modulo $2$.
\item
If $\xi^{i}\in E^{1,i}$ is an edge of $E$, then for $\xi^{i+1}:=\ol{\xi^i}$ we have that
\begin{align*}
s(\ol{\xi^i})&=r(\xi^i), \\
r(\ol{\xi^i})&=s(\xi^i).
\end{align*}
\end{itemize}
In other words, each edge $\xi$ of a bi-graph has a unique \emph{dual} $\ol\xi$ which goes in the opposite direction. A (non-ordered) pair $\{\xi^0,\xi^1\}$ with $\xi^i$ as above will be called a \emph{bi-edge}.

Given a bi-graph $E$, we shall call by \emph{underlying graph} $\E$ the non-oriented graph with the same set of axes and a edge corresponding to each bi-edge in $E$. Conversely, for each non-oriented graph $\E$ we shall call an \emph{associated bi-graph} the graph $E$ obtained by replacing each edge by a bi-edge. In what follows, we use italic characters like $E, X, Y, Z, T$ for bi-graphs and curly characters $\E, \mathcal X, \mathcal Y, \mathcal Z, \mathcal T$ for their underlying graphs and vice versa.

We shall also sometimes call a bi-graph associated to a tree a \emph{bi-tree}.

\begin{rem}\label{graph_conv_rem_edges}
The notation here is somewhat different from the one used in \cite{Iyu13}. There the labels $\xi$ and $\ol \xi$ were used as the names of the arrows rather than an operation $\xi\mapsto\ol\xi$. We employ the present notation in order to use the word ``edge'' for any kind of edges in a bi-graph. This makes it possible to choose the edges with the direction we need in a concrete part of the proof without caring whether to write the bar above or not.
\end{rem}


A bi-graph will be called a \emph{tree bi-graph} or a \emph{bi-tree} if its underlying graph is a tree graph.

If $E$ is a bi-graph and $X$ is its bi-subgraph then the edge $\xi\in E^1$ will be called \emph{adjusted to} $X$ if $\xi\notin X^1$, but $\{s(\xi),r(\xi)\}\cap X^{0}\ne\varnothing$. In other words, $\xi$ does not belong to $X$ but has at least one common axis with it. The set of all the edges in $E$ adjusted to $X$ will be denoted by $\Delta(X)$.

\subsection{Subtraction}

We introduce operations on graphs that will be extensively used throughout the paper. The first one is what we call ``subtraction'' of a bi-subgraph $X$ out of a bi-graph $E$. The na\"\i ve set-theoretical subtraction  $E$``$\setminus$''$X$, that may have been given by the data $(E^0\setminus X^0, E^1\setminus X^1, s, t)$ yields an object that in general ceases to be a graph, since the source and range maps are not well defined for it. Therefore we give the following:

\begin{defn}
Let $E$ be a bi-graph and $X$ be its bi-subgraph. Then by $E\setminus X$ we shall denote the graph in which:
\begin{itemize}
\item
$(E\setminus X)^1 = E^1\setminus X^1$, i.e.\ we remove from $E$ all the edges that belong to $X$;
\item
$(E\setminus X)^0 = \{\al\in E^0\mid \exists \xi\in (E\setminus X)^1\colon s(\xi)=\al \vee r(\xi)=\al\}$, so that we get rid of all the inner axes of $X$, i.e.\ the ones that have no adjusted edges that do not belong to $X^1$.
\item
The source and the range maps remain the same as in $E$.
\end{itemize}
\end{defn}

Observe that if $X^1=\varnothing$, i.e. $X$ is a collection of axes without edges, then in our definition $E\setminus X = E$.

\subsection{Factorization}

The second operation that we define here will be called \emph{factorization}. We will also refer to it as \emph{shrinking}.

\begin{defn}
Let $E$ be a bi-graph and let $X=\cup_{i\in I} X_i$ be its bi-subgraph, where $X_i$ are pairwise disjoint connected bi-subgraphs of $E$ and $I$ is a (countable) index set. We denote by $E/X$ the following bi-graph:
\begin{itemize}
\item
$(E/X)^0:=(E^0\setminus X^0) \cup \{\al_i\}_{i\in I}$;
\item
$(E/X)^1:=(E^1\setminus (X^1\cup\Delta(X))) \cup \Delta'(X)$, where $\Delta'(X):=\{\xi'\mid \xi\in X^1\}$ and
$$
s(\xi')=
\begin{cases}
s(\xi)& \text{if } s(\xi)\notin X^0; \\
\al_i & \text{if } s(\xi)\in X^0_i;
\end{cases}
$$
and likewise for $r(\xi)$.
\end{itemize}
\end{defn}
Informally, we ``shrink'' every connected component of a bi-subgraph $X$ to a point. In general this definition does not even preserve such a property of a graph as being locally finite, unless the connected components $X_i$ of $X$ are finite.

Observe that by the construction there is a natural injection $\iota \colon (E/X)^1 \to E^1$, given by
\begin{equation}\label{eq:iota}
\begin{cases}
\iota(\xi)=\xi & \text{ if $\xi\notin (E\setminus X)^1\cup \Delta'(X)$}; \\ \iota(\xi')=\xi & \text{ if $\xi'\in \Delta'(X)$};
\end{cases}
\end{equation}
where $\xi$, $\xi'$ are as in the definition of the factorization. By the construction we have that $\mathrm{Im}(\iota) = (E\setminus X)^1$. Hence we may consider $\iota$ as an isomorphism $\iota\colon(E/X)^1\xto{\sim} (E\setminus X)^1$.

We also denote $\tilde\iota \colon E^0 \to (E/X)^0$ the map
\[
\tilde\iota (\be) =
\begin{cases}
\be & \text{if $\be\notin X^0$};\\
\al_i & \text{if $\be\in X_i^0$}.
\end{cases}
\]
The map $\tilde{\iota}$ is obviously a surjection.

\subsection{The Valency Set}\label{subsec:Vanency}

There are several ways to define valency number for general kind of graphs. An invariant one is the following. Let $\E$ be a graph and let $\{\mathcal{X}_k\}$ be a set of finite subgraphs of $\E$ such that $\mathcal{X}_k \subseteq \mathcal X_{k+1}$ for all $k\in \Nb$ and $\bigcup_k \mathcal{X}_k = \E$. Denote by $M_k$ the set of infinite connected components of the graph $\E\setminus \mathcal{X}_k$, and for all $k$ fix an embedding $\ep_k\colon V_k\to V_{k+1}$, such that an image of a connected component $\mathcal Y_{k,\mu}$ in $V_k$ is mapped into a connected component $\mathcal Y_{k,\nu}\in V_{k+1}$ such that $\mathcal Y_{k,\nu}\subseteq \mathcal Y_{k,\mu}$. We set $V=\varinjlim_k V_k$.

Although the set $V$ depends on the choice of the sequence $\{\mathcal{X}_k\}$, any two such sets are isomorphic. Indeed, if $\{\mathcal X_k'\}$ is some other system, then there is a sequence of natural numbers $k_l$ such that $\mathcal X_{k_1} \subseteq \mathcal X_{k_2}' \subseteq \mathcal X_{k_3}\subseteq \mathcal X_{k_4}' \subseteq\dots$. Therefore there exists a sequence of embeddings $\dots\to V_{k_{4l+1}}\to V'_{k_{4l+2}} \to V_{k_{4l+3}} \to V'_{k_{4l+4}}\to \dots$ compatible with the choice of $\ep_k,\ep_k'$ that establishes an isomorphism between $V$ and $V'$. Thus we may speak about \emph{the} valency set of a graph.

We note that the valency set differs from the set of topological ends of the graph $\E$: the former is always not more than countable whereas the latter may be uncountable. There is a natural surjective map from the set of topological ends to the valency set.

Finally, we note that the valency set is stable under a factorization by a countable number of disjoint finite subgraphs. Indeed, let $\mathcal Z = \bigcup_\nu\mathcal Z_\nu$ be a graph consisting of not more than countable set of mutually disjoint finite subgraphs of $\E$. Without loss of generallity we may assume that $\mathcal X_k$ is either disjoint with $\mathcal Z_\nu$ or contains $\mathcal Z_\nu$ for all $k\in \Nb$ and all $\nu$. Set $\mathcal X'_{k}= \mathcal X_k / \bigcup_{\mathcal Z_\nu\subseteq \mathcal X_k} \mathcal Z_\nu \subset \E/\mathcal Z$ and let $V_k'$, $\ep'_k$ $V'$ be the corresponding sets of connected components, maps between them and valency set in $\E/\mathcal Z$. Then, for each $k$, there is an embedding $V'_k\subseteq V_k$ compatible with $\ep_k$. Since the graphs $\mathcal Z_k$ are disjoint, there exists an index $k'\ge k$ with an embedding $V_k\to V'_k$ compatible with the maps $\ep_l$, $\ep'_l$. By diagram chase we establish an isomorphism between $V$ and $V'$.

We denote by $\gamma(\E)$ the cardinality of the set $V$.

\section{$\K$-Theory Groups for Rose-Tree Bi-Graphs}\label{sec:Rosetree}

In this section we calculate $\K$-theory groups for a special kind of graphs, which we call rose-tree graphs. 


\begin{defn}\label{graph_defn_rose_tree}
A bi-graph $E=(E^0,E^1,s,r)$ will be called a \emph{rose-tree bi-graph} if its set of edges $E^1$ may be represented as a disjoint union $E^1=B\sqcup P$, such that:
\begin{itemize}
\item
$(E^0,B,s,r)$ is a tree bi-graph,
\item
For all $\xi \in P$ holds $s(\xi)=r(\xi)$, i.e.\ all the edges that are not in $B$ are loop edges.
\end{itemize}
The underlying graph $\E$ will be called a \emph{rose-tree} graph.
\end{defn}
In other words, a rose-tree bi-graph looks like a bi-tree which in addition may have a finite number of loop bi-edges (``rose petals'') attached to each of its axes.

Rose-tree graphs admit a simple description of their Betti and valency sets.

Indeed, for a rose-tree graph $\E$ we may establish a natural isomorphism between the set of generators of $H_1(\E)$ and the set of the ``petals'' of all the ``roses'' - the loop edges - since, by definition, form the full system of non-homotopic cycles in the graph $\E$.

As for the valency set, we may associate it with a subset of non-loop edges of $\E$ as follows. Denote by $\mathcal T$ the unique spanning tree of $\E$. Fix a root $\alpha_0$ in $\mathcal T$, and consider the orientation on $\E'$ such that the edges point out from the root. We shall say that a branch of $\mathcal T$ growing from an edge $x$ is \emph{dead} if the tree that starts at $r(x)$ and runs in the chosen orientation is finite. Let $W_0$ be the set of all edges coming out of $\al_0$ such that the branch growing from $x$ is not dead. Let $\nu_{\al_0}$ be the number of such edges. Subsequently we put an order on the axes $\alpha_k\in \mathcal T^1$ in such a way that first the first incidence level with respect to $\al_0$ is filled, then the second and so on. If there is $\nu_{\al_k}+1$ non-dead edges growing from $\al_k$, we pick up $\nu_{\al_k}$ of them and define $W_k$ to be the union of $W_{k-1}$ with this set of $\nu_{\alpha_k}$ edges. It is straightforward to check that the set $W=\bigcup_k W_k$ is isomorphic to $V$, and that $\sum_{k=0}^\infty \nu_{\al_k} = \gamma (\E)$. We shall call $\nu_{\alpha_k}$ the \emph{valency number} of $\alpha_k$.




\begin{lem}\label{graph_lem_rose_tree_K_0}
Let $E$ be a connected rose-tree bi-graph with no ``dead ends'' i.e.\ the edges with a single non-loop bi-edge attached to them. Then $K_0(\Osm_E)\cong\Zb^{\be(E)\oplus\Zb^\gm(E)}$, where $\be(E)=\dim H_1(\E)$ and $\gm(E)$ is the cardinality of the valency set of $E$.
\begin{proof}

We first observe that any bi-subgraph of a rose tree bi-graph is a disjoint union of not more than countable number of connected rose-tree bi-graphs. In particular, a connected bi-subgraph of $E$ is a connected rose-tree bi-graph.


To make calculations a bit easier and uniform, we attach a singe dead-end to the bi-graph $E$. For that we fix an axis $\al_0$ and consider the bi-graph $\tilde E$ given by
\[
\tilde E=(E^0\cup\{\al_{-1}\},E^1\cup\{x_0,\ol x_0\},\tilde r,\tilde s),
\]
where $\tilde s (\eta)=s(\eta)$, $\tilde r(\eta)=r(\eta)$ for all $\eta\in E^1$, and
\[
\tilde s(x_0)=\tilde r(\ol x_0)= \al_{-1},\quad \tilde r(x_0)=\tilde s(\ol x_0)=\al_0.
\]
By \cite[Theorem 3.1]{Iyu13}, the Cuntz-Krieger algebra associated to the bi-graph $E$ has the same $K$-theory as that of the bi-graph $\tilde E$ (Lemma \ref{graph_lem_getridoftrees_0} from the next section provides a stronger version of the cited result).

Now, for every axis $\al$ of the bi-graph $\tilde E$ except for $\al_{-1}$, we have the following list of adjusted edges:
\begin{itemize}
\item
an ``incoming'' edge $x$ with $r(x)=\al$ (resp.\ $\al = s(\ol x)$);
\item
a finite number $m_\al$ of ``outgoing'' edges $y_i$ with $s(y_i)=\al$ (resp.\ $\al = r(\ol y_i)$);
\item
a finite number $n_\al$ of bi-loops $u_j$, $\ol u_j$ with
$$
s(u_k)=r(u_k)=s(\ol u_k) = r(\ol u_k) = \al.
$$
\end{itemize}

We may write the relations occurring in the process the factorization $\Zb^{|\tilde E^1|}/(\Id-\Phi)$ by dividing the incidence matrix into blocks associated to each axis $\al$. Each of these blocks consists of the relations obtained for the edges whose target is $\al$. By an abuse of notation, we denote the generators of the group $\Zb^{|\tilde E^1|}$ by corresponding elements of $\tilde E^1$. The relation between these generators occurring at the block corresponding to the axis $\al_{-1}$ is
$$
\ol x_0 = x_0,
$$
so we have $1$ free variable.

For all the axes $\al\ne \al_{-1}$ we have the following systems of relations (for each of the edges with range $\al$):

\begin{equation}\label{graphs_eqn_rose_tree}
\left\{
\begin{aligned}
x& \colon &x &=& \sum_{j=1}^{n_\al} (u_j + \ol u_j) + \sum_{i=1}^{m_\al}y_i, \\
u_k & \colon &0 & = & \sum_{j\ne k}(u_j + \ol u_j) + \ol x + \sum_{i=1}^{m_\al}y_i & &\text{($n_\al$ relations}),\\
\ol u_k & \colon& 0 & = & \sum_{j\ne k}(u_j + \ol u_j) + \ol x + \sum_{i=1}^{m_\al}y_i & &\text{($n_\al$ relations}),\\
\ol y_l& \colon &\ol y_l & =&  \ol x + \sum_{j=1}^{n_\al} (u_j + \ol u_j) + \sum_{i\ne l}y_i & & \text{($m_\al$ relations).}
\end{aligned}
\right.
\end{equation}
Observe that the relations obtained for $u_j$ and $\ol u_j$ coincide, so that in fact the system of relations for the edges with range $\al$ is equivalent to:
\[
\left\{
\begin{aligned}
x &=& \sum_{j=1}^{n_\al} (u_j + \ol u_j) + \sum_{i=1}^{m_\al}y_i, \\
0 & = & \sum_{j\ne k}(u_j + \ol u_j) + \ol x + \sum_{i=1}^{m_\al}y_i, \\
\ol y_l & =&  \ol x + \sum_{j=1}^{n_\al} (u_j + \ol u_j) + \sum_{i\ne l}y_i.
\end{aligned}
\right.
\]
Subtracting the first relation from all the others, we rewrite this system as follows:
\[
\left\{
\begin{aligned}
x &= \sum_{j=1}^{n_\al} (u_j + \ol u_j) + \sum_{i=1}^{m_\al}y_i, \\
-x & =   -(u_k+\ol u_k) + \ol x, \\
-x+ \ol y_l & =  \ol x +  - y_l;
\end{aligned}
\right.
\]
or, in a simplified form:
\[
\left\{
\begin{aligned}
x &= \sum_{j=1}^{n_\al} (u_j + \ol u_j) + \sum_{i=1}^{m_\al}y_i, \\
(u_k+\ol u_k) & =   x + \ol x,  \\
y_l + \ol y_l & = x+ \ol x.
\end{aligned}
\right.
\]
Let us now calculate the number of free variables that appear at the axis $\al_0$ with the "incoming" edge $x_0$. We have already established that $\ol x_0 = 0$, so we obtain $2m_{\al_0} + 2n_{\al_0} + 1$ variables and $m_{\al_0} + n_{\al_0}$ relations. Since all the coefficients equal $1$, this gives us in total $m_{\al_0} + n_{\al_0}$ free variables.

For all the subsequent steps we may assume the values of $x$ and $\ol x$ to be already fixed (these $x$'s were some of the $y$'s form the previous step). Therefore at each axis $\al$ we shall have $2m_\al +2n_\al$ new variables and $m_\al + n_\al + 1$ relations, resulting in $m_\al + n_\al - 1$ additional free variables at each axis.

Thus, every axis gives the following number of free variables to the pool of generators of the Grothendieck group:
\begin{itemize}
\item $m_\al$ of those corresponding to the loop bi-edges; the set of all these generators is naturally isomorphic to the set of bi-loops of $E$;
\item $1$ for the axis $\al_{-1}$ and $n_\al-1$ for all $\al\ne \al_{-1}$ corresponding to non-loop bi-edges; observe that these numbers are precisely the valency numbers $\nu_\al$ of the axes we have defined above.
\end{itemize}

The total set of generators thus equals $\beta(\E)+\gamma(\E)$ so that $\K(\Osm_\E)=\Zb^{\beta(\E)+\gm(\E)}$, as it has been claimed.


\end{proof}
\end{lem}

In this proof, both the assumption for $E$ to have no ``dead ends'' as well as the adjoining of an additional root axis with a bi-edge were only made to make the calculations look more uniform so that unnecessary technicalities may be avoided. In the next section we shall see that the requirement for a rose-tree graph to have no ``dead ends'' may be relaxed.

The calculation of the Whitehead groups for the algebras associated to rose-tree bi-graphs is also direct. For this we do not even need to restrict ourselves to bi-graphs with no ``dead ends''. Namely:

\begin{lem}\label{graph_lem_rose_tree_K_1}
Let $E$ be a rose-tree bi-graph. Then $K_1(\Osm_\E)\cong\Zb^{\be(\E)}$.
\begin{proof}


Recall that $K_1(\Osm_\E)\cong\ker(\Id-\Phi)$. Suppose that $\xi\in \ker(\Id-\Phi)$, $\xi=\sum_{e\in E_\xi} k_x x$ where $E_\xi\subset E^1$ is a finite subset of edges and $k_e\in \Zb\setminus\{0\}$. Our goal is to show that $\xi$ necessarily has the form $\xi = \sum k_z(z-\ol z)$ where $z$ are loop edges of the bi-graph $E$.



To prove this, we fix a root axis in $E$, and proceed by induction in the axes of $E_\xi$, starting with those that lay furthest from the root. More precisely, let first $\al$ be an axis of $E_\xi$, such that there is no other edge $\al'$ that lies further from the root and no such edge $e$ which comes as a term with a non-zero coefficient in $\xi$, such that $r(e)=\al$, $s(e)=\al'$ or vice versa. Then we may write
\[
\xi = k_x x + k_{\ol x}\ol x + \sum_{s(z)=r(z)=\al}(k_z z + k_{\ol z}\ol z) + \eta_0,
\]
where $x$, $r(x)=\al$ is the only ``outgoing'' non-loop edge, $\ol x$ is its opposite, $z$ are ``real'' loop edges with $\ol z$ being their opposites, and $\eta_0$ consists of the rest terms. Let $y_1,\dots,y_l\ne \ol x$ be those edges of $E$ for which $r(y_j)= \al$.

By definition, we have that
\[
\begin{aligned}
(\Id-\Phi) x & = x + \ol x - \sum_{s(z)=r(z)= \al}(z+\ol z) - \sum_{j=1}^ly_j, \\
(\Id-\Phi)z & = -\sum_{z'\ne z} (z'+\ol z') - \sum_{j=1}^l y_j, \\
(\Id-\Phi)\ol z & = -\sum_{z'\ne z} (z'+\ol z') -\sum_{j=1}^l y_j.
\end{aligned}
\]
The edge $x$ and the loop edges with source and range $\al$ are the only ones among the edges $e\in E^1$ with nonzero coefficients $k_e\ne 0$ in $\xi$ and for which $(\Id-\Phi)e$ contains some of the terms for the edges $y_j$ with nonzero coefficients. Thus, we have that
\[
(\Id-\Phi)\xi = -\left(k_{x}+ \sum_{s(z)=r(z)=\al}(k_z+k_{\ol z})\right)\sum_{j=1}^l y_j + \eta_1,
\]
where $\eta_1$ is formed by the rest terms. Since $\xi\in \ker(\Id-\Phi)$ we obtain
\begin{equation}\label{eq:K_1:k_x}
k_{x} = \sum_{s(z)=r(z)=\al}(-k_z-k_{\ol z}).
\end{equation}
From the other hand
\[
(\Id-\Phi)\xi = \left(\sum_{z'\ne z}(-k_{z'}-k_{z'})- k_{x} \right)(z+\ol z) + \eta_2,
\]
where $\eta_2$ is formed by the terms that do not include any loop edges attached to $\al$. Again, since $\xi\in \ker(\Id-\Phi)$, we get
\begin{equation}\label{eq:K_1:k_x:k_z}
\sum_{z'\ne z}(-k_{z'}-k_{z'})- k_{x} = 0,
\end{equation}
for all $z$ as above. Summing this equality for given $z$ with \eqref{eq:K_1:k_x}, we get simple relations
\begin{equation}\label{eq:K_1:k_z+olk_z}
k_z+k_{\ol z} =0,
\end{equation}
for all loops $z$. Substituting \eqref{eq:K_1:k_z+olk_z} into \eqref{eq:K_1:k_x:k_z}, we see that $k_{x}=0$. Finally
\[
(\Id-\Phi)\xi = \left(k_{\ol x}+ \sum_{s(z)=r(z)=\al} (k_z+k_{\ol z})\right) \ol x + \eta_3,
\]
so that
\[
k_{\ol x}=-\sum_{s(z)=r(z)=\al} (k_z+k_{\ol z}) = 0.
\]
Hence, for the ``utmost'' edges from the root, the only terms with nonzero coefficients are to be of the form $z$ or $\ol z$, with $r(z)=s(z)=\al$, and it holds that $k_z=-k_{\ol z}$.

Let now $\al$ be such an axis that there are edges in $E_\xi$ with $\al$ being their source or range, and that for all axes $\al'$ lying further from the root than $\al$ the relations from the above paragraph hold. As above, we denote by $y_1,\dots,y_l$ the non-loop edges with $r(y_j)=\al$ and $s(y_j)$ lying further from the root than $\al$. We also set $x$ with $s(x)=\al$ to be the edge pointing towards the root and $\ol x$ be its opposite. By the above hypothesis, we have that $k_{y_j}=0$ and $k_{\ol y_j}=0$ for all $j=1,\dots,k$. Thus, we again have that
\[
(\Id-\Phi)\xi = \left(k_x+\sum_{s(z)=r(z)=\al}(k_z+k_{\ol z}) \right)\sum_{j=1}^l y_j + \eta_4,
\]
so that
\[
k_x = \sum_{s(z)=r(z)=\al}(-k_z-k_{\ol z}),
\]
and, from the other hand,
\[
(\Id-\Phi)\xi = \left(\sum_{z'\ne z}(-k_{z'}-k_{z'})- k_{x} \right)(z+\ol z) + \eta_5.
\]
Hence we again obtain that $k_z=-k_{\ol z}$ and $k_x=0$. Finally,
\[
(\Id-\Phi)\xi = \left(k_x+ k_{\ol x} - \sum_{s(z)=r(z)=\al}(k_z+k_{\ol z}\right) - \sum_{j=1}^l k_{\ol y_j})\ol x + \eta_6,
\]
and since $k_{\ol y_j}=0$ for all $1\le j \le l$ by the hypothesis we have that $k_{\ol x}=0$ in order for $\xi$ to be in $\ker(\Id-\Phi)$. Thus, there are no non-loop edges $e$ with $\al$ being its source or target, such that $k_e\ne 0$, and for all $z$ with $r(z)=s(z)=\al$ we have $k_z=-k_{\ol z}$. Proceeding by induction, we prove it for all the axes $\al\in E^0$, such that $\al$ is a source or a target of an edge of $e\in E_\xi$. Clearly, the proof does not depend on the choice of the root we've made as the first step.

Summarizing, we obtain that if $\xi\in \ker(\Id-\Phi)$, then $\xi$ may be represented as a finite sum $\xi = \sum k_z(z-\ol z)$, where $s(z)=r(z)$. The group $\mathrm{K}_1(\Osm_\E)$ is a free abelian group whose generators stand in 1-to-1 correspondence with the loop bi-edges of the graph $E$, i.e. $\K_1(\Osm_\E)=\Zb^{\be(\E)}$. QED.
\end{proof}
\end{lem}

\section{Technical Lemmata}\label{sec:TechLemmata}

The two lemmata we shall now prove are in fact results from linear algebra that may be useful on their own. Their main purpose in the current article is to prove the preservation of $K_0$ and $K_1$ of the associated algebra under ``shrinking'' of a countable number of disjoint finite trees in a locally finite graph.

\begin{lem}\label{graph_lem_getridoftrees_0}
Let $G$ and $H$ be abelian groups and $T\colon G\oplus H \to G\oplus H$ be an endomorphism. Set $P\colon G\oplus H \to G\oplus H$ to be a map for which $P|_G = \Id_G$ and $P(x) = T(x)$ for any $x\in H$. Suppose also that for any $x\in H$ there exists a finite number $n(x)\in\Nb$ such that $P^{n(x)}(x)\in G$ and define $P^\iy(x) := P^{n(x)}(x)$. Finally, let $\tilde T\colon G\to G$ be an endomorphism of $G$ defined as $\tilde T := P^\iy \circ T$. Then there is a natural isomorphism
\begin{equation}\label{eq:J_coker}
J\colon G/\Im (\Id_G -\tilde T) \xto[\sim]{} (G\oplus H)/\Im (\Id_{G\oplus H} - T),
\end{equation}
given by
\begin{equation}\label{eq:getridoftrees_1_J}
J(u + \tilde T(G)) = u + T(G\oplus H).
\end{equation}
\begin{proof}


\textbf{1. The map $J$ is well-defined.} To show this, we need to prove that the image of the $0$-class in $G/\Im (\Id_G -\tilde T)$ is the $0$-class in $(G\oplus H)/\Im (\Id_{G\oplus H} -T)$.

Indeed, let $p$ be a natural projection $p\colon G\to G/\Im (\Id_G -\tilde T)$. Then $p(u)=0$ iff there exists such $v\in G$ that $u= v - P^\iy T(v)$. Suppose that $T(v)= g_0+h_0$, where $g_0\in G_0$ and $h_0\in H$. Observe that by the construction
$$
P^\iy T(v) = P^\iy (g_0 + h_0) = g_0 + P^\iy(h_0).
$$
Also, by the conditions of the lemma, we have that $P(h_0) = T(h_0)$ and $P^\iy(h_0)=P^{n(h_0)}$. Thus,
$$
P^\iy T(v) = g_0 + P^{n(h_0)-1}T(h_0).
$$
Now let $T(h_0)= g_1 + h_1$ with $g_1\in G$ and $h_1\in H$. Then we may analogously rewrite
$$
P^\iy T(v) = g_0 + g_1 + P^{n(h_0)-2}T(h_1).
$$
Continuing this process, we obtain
$$
P^\iy T(v) = g_0 + g_1 + \dots + g_{n(h_0)-1} + g_{n(h_0)},
$$
where, by the construction, $g_j = T(h_{j-1}) - h_j$ and $g_{n(h_0)} = T(h_{n(h_0)-1})$ are elements of $G$.

Now, set
$$
w:= v+ h_0 + h_1 + \dots + h_{n(h_0)-1}.
$$
We have that
\[
\begin{split}
w-T(w) =& v+ h_0 + h_1 + \dots + h_{n(h_0)-1} - T(v+ h_0 + h_1 + \dots + h_{n(h_0)-1}) = \\
& v - T(v) + h_0 - T(h_0) + h_1 - T(h_1) + \dots + h_{n(h_0)-1} - T(h_{n(h_0)-1}) = \\
& v - g_0 - g_1 - \dots - g_{{n(h_0)-1}} - g_{n(h_0)} =\\
& v - P^\iy T(v) = u.
\end{split}
\]
Therefore $u\in\Im(\Id_{G\oplus H} - T)$. Thus $J(\Im (\Id_G -\tilde T)) \subseteq \Im(\Id_{G\oplus H} - T)$, and so the map $J$ is well-defined.

\textbf{2. The map $J$ is injective.} Indeed, let $u\in G$ be such an element that $J(u+\Im(\Id_G - \tilde T)) = \Im(\Id_{G\oplus H} - T)$. By definition of $J$ this means that $u\in \Im(\Id_{G\oplus H} - T)$, or, in other words, $u=v-T(v)$ for some $v\in G\oplus H$. The element $v$ admits a unique representation $v=g+h$, where $g\in G$, $h\in H$. Observe that
$$
T(g) = T(g+h-h) = T(v) - T(h) = v-u - T(h) = g-u + h - T(h).
$$
Thus
\[
\begin{split}
g - \tilde T(g) = & g - P^\iy(g-u + h - T(h)) =\\
&  g-g+u + P^\iy(h-T(h))=\\
&  u + P^\iy(h) - P^\iy(h)=\\
&  u.
\end{split}
\]
Therefore $u\in\Im(\Id_G - \tilde T)$, and so the map $J$ is injective.

\textbf{3. The map $J$ is surjective}

Indeed, let $v = g + h\in G\oplus H$. Consider the element $P^\iy(v)$. By the conditions, we have that $P^\iy(v)\in G$. We argue that $v$ and $P^\iy(v)$ define the same element in $G/\Im(\Id_{G\oplus H} -  T)$, or, in other words, $v-P^\iy v \in \Im(\Id_{G\oplus H} -  T)$.

Indeed, represent $v=g+h$ as above. Then
$$
v-P^\iy(v) = g+h - P^\iy(g+h) = g+h -g - P^\iy(h) = h-P^\iy(h).
$$
Observe that by definition of $P$ there exists a number $n$ such that $P^\iy(h) = P^n(h)$. Therefore
$$
h - P^\iy(h) = h - P^n(h) = h - P(h) + P(h) - P^2(h) + \dots + P^{n-1}(h) - P^n(h).
$$
Write $P^j(h)=:g_j + h_j$ with $g_j\in G$ and $h_j\in H$. Then
$$
P^j(h) - P^{j+1}(h) = g_j + h_j - P(g_j + h_j) = h_j - P(h_j) = h_j - T(h_j),
$$
since by the construction $P(h_j) = T(h_j)$.
Hence we may rewrite
\[
\begin{split}
v - P^\iy(v) =& h - P^\iy(h) =\\
& h - T(h) + h_1 - T(h_1) + \dots + h_{n-1} - T(h_{n-1}) =\\
&  (h +\sum_{j=1}^{n-1}h_j) - T (h+\sum_{j=1}^{n-1}h_j) \in \Im(\Id_{G\oplus H} - T).
\end{split}
\]
Thus
$$
v+ \Im(\Id_{G\oplus H} - T) = P^\iy(v) + \Im(\Id_{G\oplus H} - T) = J(P^\iy(v) +  \Im(\Id_G - \tilde T)),
$$
so that the map $J$ is surjective.\newline

Putting (1), (2) and (3) together, we obtain that $J$ is an isomorphism.
\end{proof}
\end{lem}

The next lemma establishes an analogous result for the kernel.

\begin{lem}\label{graph_lem_getridoftrees_1}
In the conditions of Lemma \ref{graph_lem_getridoftrees_0}, there is an isomorphism \[\ker (\Id_{G\oplus H} - T) \xto[\sim]{} \ker (\Id_G - \tilde T),\]
given by the restriction of the natural projection:
\begin{equation}
\Pi\colon G \oplus H \to G.
\end{equation}
\begin{proof}
We want to show that $\Pi$ restricts to an isomorphism $\ker (\Id_{G\oplus H} - T)\to \ker (\Id_G - \tilde T)$, which we shall abusively denote by $\Pi$.

\textbf{1. The restriction of $\Pi$ is well-defined.}

Let $g+h\in \ker(\Id_{G\oplus H} - T)$. This means that
$$
g+h-T(g)-T(h)=0.
$$
Applying $P^\iy$ to both sides of the equation we obtain:
\[
\begin{split}
0 =& P^\iy(g+h-T(g)-T(h)) = \\
&  g - P^\iy(h) - \tilde T(g) - P^\iy\circ P(h) = \\
&  g- \tilde T(g),
\end{split}
\]
and so we indeed have that $g\in \ker(\Id_G -\tilde T(g))$.

\textbf{2. The restriction of $\Pi$ is surjective.}

Let $g\in \ker(\Id_G -\tilde T)$, so that $g-\tilde T(g) = 0$. We need to show that there exists such $h$ that $g+h - T(g+h) =0$.

To see this, recall that $\tilde T(g) = P^\iy T (g) = P^n T(g)$ for some $n$, so that we may write
\[
\begin{split}
0 & = g - \tilde T(g) \\
& = g - P^n T (g) \\
&= g - T(g) + T(g) - P\cdot T(g) + P\cdot T(g) -\dots - P^{n-1}\circ T(g) + P^{n-1}\circ T(g) - P^n\circ T(g).
\end{split}
\]
Let $T(g) = g_0 + h_0$, and denote $g_k + h_k := P(h_{k-1})$, where $g_k\in G$ and $h_k\in H$ for all $k = 0,\dots, n-1$. Then
\[
\begin{split}
P^{k-1} \circ T(g) - P^{k} \circ T(g)  = &
P^{k-1} (g_0 + h_0 - P (g_0 + h_0)) = \\
& P^{k-1} (h_0 - P(h_0)) = \\
& P^{k-1} (h_0 - g_1 - h_1) = \\
& -g_1 + P^{k-2}(h_1 + g_1 - h_2 - g_2) = \\
& P^{k-2}(h_1 - g_2 - h_2) = \\
& \dots \\
& P (h_{k-1} - g_k - h_k) = \\
& -g_k + g_k + h_k - P(h_k) = \\
& h_k - P(h_k).
\end{split}
\]
Hence we have
\[
\begin{split}
0  = & g - T(g) + T(g) - P\cdot T(g) + P\cdot T(g) -\dots - P^{n-1}\cdot T(g) + P^{n-1}\cdot T(g) - P^n \cdot T(g) = \\
& g - T(g) + h_0 - P(h_0) + h_1 - P(h_1) + \dots + h_n - P(h_n) = \\
& g - T(g) + (h_0+h_1+\dots +h_n) - T(h_0 + h_1 + \dots +h_n).
\end{split}
\]
Set $h := h_0 + h_1 + \dots h_n$. Then
$$
g + h \in \ker (\Id_{G\oplus H} - T),
$$
and $\Pi(g+h) = g$. Hence the map $\Pi$ is surjective.

\textbf{3. The restriction of $\Pi$ is injective.}

Indeed, let $g\in \ker(\Id_G -\tilde T(g))$ have two preimages $g+h',g+h''\in \ker(\Id_{G\oplus H} - T)$. Then $h:=h'-h''\in \ker(\Id_{G\oplus H} - T)$, so that
$$
h=T(h)=P(h).
$$
where the second equality holds since $h\in H$. Applying $P$ to both sides of the equation $h=Ph$, we obtain
$$
P(h)=P^2(h),
$$
and thus
$$
h = P^2(h).
$$
Analogously $h=P^n(h)$ for all $n\in \Nb$, and so
$$
h=P^\iy(h).
$$
But $P^\iy(H)\in G$, hence $h\in G\cap H = \{0\}$, and therefore $h=0$. Thus $\Pi$ is injective.

Putting all the above results together, we see that the restriction of the map $\Pi$ is indeed an isomorphism. QED.
\end{proof}
\end{lem}

The above lemmata yield the following

\begin{prop}\label{graph_cor_getridoftrees}
Let $E$ be a locally finite bi-graph and let $X = \bigcup_{i\in I} X_j$ be a bi-subgraph of $E$, where $X_i$ are mutually disjoint finite bi-tree subgraphs of $E$. Then
$$K_i(\Osm_\E) = K_i(\Osm_{\E/\mathcal X}),\quad i=0,1.$$
\begin{proof}

We abusively denote by $\xi$ the generators of $\Zb^{|E^1|}$ corresponding to the edges $\xi\in E^{1}$. Set $P\colon \mathbb{Z}^{|E^1|} \to \mathbb{Z}^{|E^1|}$ to be such a map that $P(\xi) = \Phi(\xi)$ for $\xi\in X^1$ and $P(\xi) = \xi$ otherwise.

Now, since $X_i\cap X_j = \varnothing$ for $i\neq j$, we have that for all $\et \in X$ there is a unique index $i$ for which $\et\in X_i$. Let $\al := r(\et)$. Then we have that
$$
P(\et) = \Phi(\et) = \sum_{s(\xi) =\al}\xi + \sum_{s(\et)=\al} \et,
$$
where $\xi_\al\in E^1\setminus X^1$ and $ \et\in X_j^1$. Observe that since
\begin{itemize}
\item
the expression for $\Phi(\et)$ does not contain $\ol\et$;
\item
the double graph $X_j$ is a bi-tree, therefore it contains no loops other than $(\{s(\et),r(\et)\},\{\et,\ol\et\},r,s)$;
\end{itemize}
we shall encounter no such $\et$ in the expression for $P^n(\et)$ for any $n\in\mathbb{N}$. Since $ X_j$ is finite, this means that there exists such $n\in\mathbb{N}$ that $P^n(\xi)\in \Zb^{|E^1\setminus  X^1|}$.

Thus, the groups $\mathbb{Z}^{|E^1\setminus  X^1|}=:G$, $\mathbb{Z}^{| X^1|}=:H$ and the map $\Phi\colon G\oplus H \to G\oplus H$ satisfy the conditions of the Lemma \ref{graph_lem_getridoftrees_0}. QED.
\end{proof}
\end{prop}

Informally, this result means that, given two locally finite graphs $\E$ and $\tilde{\E}$, such that $\tilde \E$ is obtained from $\E$ by ``shrinking'' a countable number of finite non-intersecting bi-tree subgraphs, the $\K$-theory groups of the algebras $\Osm_{\E}$ and $\Osm_{\tilde \E}$ will be naturally isomorphic.

\begin{rem}
Observe that, as we have announced earlier, Proposition \ref{graph_cor_getridoftrees} makes it possible for us to relax the condition for the rose-tree to have no dead ends. Indeed, for each axis $\al$ of a bi-graph $E$, consider the maximal finite bi-subgraph $X_\al$ such that $E\setminus X_\al$ is connected. This bi-subgraph $X_\al$ is evidently a rose-tree. Let $Y_\al$ be the spanning bi-tree of $X_\al$. Obviously, for different axes $\al$ and $\al'$ the graphs $Y_\al$ and $Y_\al'$ do not intersect, and so the algebra $\Osm_{E/(\bigcup_{\al\in E^0}Y_\al)}$ has the same $\mathrm K$-theory as $\Osm_{E}$. By construction, $E/(\bigcup_{\al\in E^0}Y_\al)$ is a rose-tree bi-graph, and the reader may check, that it has no dead ends.
\end{rem}

\section{Reduction Theorem}\label{sec:Reduction}

In this section we are proving the main result of the paper. The goal of the theorem is to show that for an arbitrary locally finite infinite bi-graph $E$ there exists a rose-tree graph $\tilde{E}$, such that the $\K$-groups of their associated algebras are naturally isomorphic. Moreover, the bi-graph $\tilde E$ will be constructed out of $E$ by means of factorization by finite bi-tree graphs in several steps. That will allow us to establish a connection of $\K_i(\Osm_\E)$ with first homology group of $\E$ and its valency set. This is done in the next section.

\begin{thm}\label{thm:reduction}
Let $E$ be an infinite locally finite bi-graph. Then there is a rose-tree bi-graph $\tilde E$ such that there is are natural isomorphisms $\K_i(\Osm_\E)\cong \K_i(\Osm_{\tilde \E})$ for $i=0,1$.


\begin{proof}

Our goal is to construct a sequence of bi-graphs starting with $E$, such that one may be obtained from another by factorizing by a countable number of non-intersecting bi-tree subgraphs, with the last one being a rose-tree bi-graph with no dead ends. According to the Proposition \ref{graph_cor_getridoftrees}, all these graphs will have isomorphic $\mathrm K$-theory groups. 
We split the proof into several steps.

\textbf{Step 1:}

First we take a connected finite bi-subgraph $X_0$ of $E$. Then we choose a finite bi-subgraph $X_1$ of $E$ with the following properties:
\begin{itemize}
\item
$X_0^1\cap X_1^1 = \varnothing$.
\item
All the edges adjusted to $X_0$ lie in $X_1$. More formally, if $\xi\in E^1\setminus X^0_1$ is such that $s(\xi)\in X_0^0$ or $r(\xi)\in X_0^0$, then $\xi\in X_1^1$. We denote the graph composed out of these double edges and their sources and ranges by $\Delta_1$.
\item
The bi-graph $X_1$ has only one connected component in each connected component of $E\setminus X_0$. In other words, if there are two or more connected components of the graph $\Delta_1$ in a single connected component of $E\setminus X_0$, then there is a single connected component of $X_1$ containing all these components of $\Delta_1$.
\end{itemize}

Observe that by the construction the graphs $X_0$ and $X_1$ intersect only by axes, and moreover $X_0\cap (E\setminus (X_0\cup X_1)) = \varnothing$.

Now take the graph $X_0\cup X_1$ and choose a finite subgraph $X_2$ in $E$ the same way we have chosen $X_1$ for the bi-subgraph $X_0$. Proceeding inductively, we set $X_{n+1}$ to be a bi-subgraph of $E$ satisfying the three above conditions for the bi-subgraph $\bigcup_{i=1}^n X_i$.

The main property of this representation is that:
\begin{itemize}
\item
$$X_i \cap X_{i+2} = \varnothing.$$
Indeed, by the construction, all the edges adjusted to $X_i$ belong either to $X_{i-1}$ or $\Delta_{i+1}\subseteq X_{i+1}$.
\item
$$
\bigcup_{i=1}^{\iy}X_i = E.
$$
Indeed, pick any axis $\al$ in $X_0$ and let $\Gamma_n$ be the $n$'th incidence bi-subgraph of $E$ with respect to $\al$. Then, since $\E$ is connected, we have by the construction that  $E=\bigcup_n \Gamma_n \subseteq X_n \subseteq E$.
\end{itemize}
Now, after having the graph $E$ represented as such a union of finite subgraphs, we may start the reduction of $E$ to a rose-tree bi-graph.

\textbf{Step 2:}

For each connected component $X_{2k+1,\mu}$ of the graph $X_{2k+1}$ (i.e.\ of those with \emph{odd} indices) we choose its spanning bi-tree $T_{2k+1,\mu}$. 
The bi-trees $T_{2k+1,\mu}$ are finite and pairwise disjoint for all $k\in \Nb$ and $\mu$ as the graphs $X_{2k+1}\cap X_{2k+3}=\varnothing$ for all $k\in \Nb \cup \{0\}$, therefore we are in the conditions of Proposition \ref{graph_cor_getridoftrees}. Thus we have that $K_i(\Osm_\E)=K_i(\Osm_{\hat \E})$, where $\hat E = E/\bigcup_{k,\mu} T_{2k+1,\mu}$ and $i=0,1$.

The graph $\hat E$ consists of the ``remnants'' $\hat X_{k,\mu}$ of the graphs $X_{k,\mu}$ of two sorts. First are the ``remnants'' $\hat X_{2k+1,\mu}$ of the graphs $X_{2k+1,\mu}$ in $\hat E$ such that
\begin{itemize}
\item
Each graph $\hat X_{2k+1,\mu}$ has a single axis that we denote by $\al_{2k+1,\mu}$.
\item
The graph $\hat X_{2k+1,\mu}$ has $\be(\mathcal X_{2k+1,\mu})$ bi-edges, and all of them are obviously loops. Here $\mathcal X_{2k+1,\mu}$ is the underlying graph of $X_{2k+1,\mu}$.
\end{itemize}

The graphs $\hat X_{2k,\nu}$ are obtained from the graphs $X_{2k, \nu}$ in the following way: if $\xi\in X_{2k,\nu}^1$ is such that $s(\xi)\in X_{2k\pm 1, \la}^0$ (resp. $r(\xi)\in X_{2k\pm 1,\mu}^0)$, then the edge $\xi$ is replaced by an edge $\hat \xi$, such that $s(\hat \xi)= \al_{2k\pm 1,\la}$ (resp. $r(\hat\xi) = \al_{2k\pm 1,\mu}$). In other words, for each connected component of $X_{2k\pm 1,\mu}$ such that $X_{2k, \nu}\cap X_{2k\pm 1,\mu} \ne \varnothing$ all the points of $X_{2k, \nu}\cap X_{2k\pm 1,\mu}$ are "glued together" in $\al_{2k \pm 1,\mu}$. An analogous reassignment is performed for the dual edges $\ol \xi$.

\textbf{Step 3:}

We construct a graph $\check E$ out of $\hat E$ by ``splitting'' the axes $\al_{2k+1,\mu}$ in two and connecting them with a bi-edge. More precisely we replace the axes $\al_{2k+1,\mu}$ belonging to the graphs $\hat X_{2k,\nu}$ by the axes that we denote by $\al_{2k,\mu}$, and add a double edge $\{\et_{2k+1,\mu},\ol \et_{2k+1,\mu}\}$ such that $s(\et_{2k+1,\mu})=\al_{2k,\mu}$ and $r(\et_{2k+1,\mu})=\al_{2k+1,\mu}$. We denote these graphs with replaced axes by $\check X_{2k,\nu}$.

The graph $\check E$ has the same $\K$-theory as $\hat E$. Indeed, the bi-graphs
$$
\check Y_{2k+1,\mu} :=(\{\al_{2k+1,\mu},\al_{2k,\mu}\}, \{\et_{2k+1,\mu}, \ol \et_{2k+1,\mu}\},s,r)
$$
are non-intersecting tree bi-graphs (each consisting of a single bi-edge), and we obviously have a countable number of them. Therefore we may apply Proposition \ref{graph_cor_getridoftrees} to $\check E$ and $\bigcup_{2k+1,\mu} Y_{2k+1,\mu}$. But now $\hat E = \check E/(\bigcup_{2k+1,\mu}Y_{2k+1,\mu})$, and so $\K_i(\Osm_{\hat \E}) \cong \K_i(\Osm_{\check \E})$, $i=0,1$.

\textbf{Step 4:}

In the graph $\check E$, denote
$$
\check Z_{2k+1,\mu}:=\hat X_{2k+1,\mu} \cup \bigcup\{\hat X_{2k+2,\nu}\mid \al_{2k+1,\mu}\in \hat X_{2k+2,\nu}^0\}
$$
In each of the bi-graphs $\check Z_{2k+1,\mu}$ choose a spanning tree bi-graph $\check T_{2k+1,\mu}$. 
We also choose a spanning bi-tree $\check T_0$ in $\hat X_0$.

The graphs $\check T_{2k+1,\mu}$ and $\check T_0$ are non-intersecting finite bi-trees. Otherwise we would have had that $\al_{2k+3,\mu}\in \hat X_{2k+1, \nu}$ for some indices $\mu,\nu$, but we have replaced all the edges of the form $\al_{2k+3,\mu}$ in $\hat X_{2k+1, \la}$ by $\al_{2k+2,\mu}$ at the previous step, so this is not possible. Because of that, we may use Proposition \ref{graph_cor_getridoftrees} once again, and obtain that, for the graph
$$
\tilde E := \check E \;/ \bigcup_{2k+1,\mu} T_{2k+1,\mu},
$$
we have that
$$
\K_i (\Osm_{\tilde \E}) \cong K_i(\Osm_{\check \E}).
$$
But, by the construction, we know that
$$
\K_i(\Osm_{\check \E}) \cong \K_i(\Osm_{\hat \E}) \cong \K_i(\Osm_\E),
$$
so that
$$
\K_i (\Osm_{\tilde \E})\cong K_i (\Osm_\E).
$$
Here $i=0,1$.

\textbf{Step 5:}

It remains to show that the bi-graph $\tilde E$ is a rose-tree bi-graph. Observe that the graphs $\check Z_{2k+1,\mu}$ and $\hat X_0$ are replaced by the graphs $\tilde Z_{2k+1,\mu}$ and $\tilde X_0$ respectively, each of which contains only a single axis which we shall identify with $\al_{2k+1,\mu}$ for $\tilde Z_{2k+1,\mu}$ and denote by $\al_0$ for $\tilde X_0$.

The bi-graphs $\check Y_{2k+1,\mu}$ become the bi-graphs
$$
\tilde Y_{2k+1,\mu} = (\{\al_{2k-1,\mu}, \al_{2k+1,\mu}\}, \{\xi_{2k+1,\mu},\ol \xi_{2k+1,\mu}\},s,r).
$$
Here the index $(2k-1,\nu)$ is the same as the index of the component $\check Z_{2k-1,\nu}$ to which the axis $\be_{2k+1,\mu}$ of the graph $\check Y_{2k+1,\mu}$ belongs. For $k=0$ this is the axis $\al_0$.

By the construction, the bi-edges $\{\et_{2k+1,\mu},\ol \et_{2k+1,\mu}\}$ of the graphs $\tilde Y_{2k+1,\mu}$ are the only non-loop bi-edges of the bi-graph $\tilde E$. It also follows from our construction that there are no other simple cycles in the graph $\tilde \E$ other than loops. Indeed, one may check that if there were any then there should have existed two disjoint bi-subgraphs of the form $X_{k,\mu}$, $X_{k,\nu}$, lying in the same connected component of $E\setminus X_{k-1}$, which contradicts of the construction, see Step 1. 


Therefore $\tilde E$ is a rose-tree bi-graph, and by Lemmata \ref{graph_lem_rose_tree_K_0} and \ref{graph_lem_rose_tree_K_1} there is are natural isomorphisms
\[
\K_0(\Osm_\E)\cong \Zb^{\beta(\tilde \E )}\oplus \Zb^{\gamma(\tilde \E)},\quad \K_1(\Osm_\E)\cong \Zb^{\beta(\tilde E)}.
\]
\end{proof}
\end{thm}

\section{Connection with topological invariants of the underlying graph}\label{sec:Cnnection}

We shall finally show that there is an invariant description of the generators of $\K_i(\Osm_\E)$ in terms of topological invariants of the graph $\E$ itself. The claims look almost obvious, yet their proofs require a substantial amount of technicalities arising from the fact that we are need to work with free groups having a countable number of generators. Thus, in order to show the naturality of isomorphisms, we may not simply count the number of generators, but are to refer to particular objects in the graph. These results are far from being elegant, their only purpose is to keep track where do the generators of $\K$-groups come from.

We start with the description of the Whitehead group as a simpler task.

\begin{cor}\label{cor:inv:K_1}
There is a natural isomorphism $\K_1(\Osm_\E)\cong H_1(\E)$, where $H_1(\E)$ is the first homology group of a realization of $E$ considered as a simplicial complex.
\begin{proof}

Let first $E$ be some bi-graph and $X$ be a collection of non-intersecting rose-tree bi-subgraphs of $E$. Suppose now that $y_j$, $j=0,\dots, m$ is collection of edges in $E/X$ such that $s(y_m)= r(y_0)$ and $s(y_j)=r(y_{j+1})$ for the rest. Then one may check directly that $\sum_{j=1}^m (y_j-\ol y_j)\in \ker(\Id_{E/X}-\Phi_{E/X})$. Moreover, \emph{the} preimage of this sum in $\Zb^{|E^1|}$ under the map $\Pi$ as in Lemma \ref{graph_lem_getridoftrees_1} is of the form $\sum_{j=0}^n(x_k-\ol x_k)$ where again $r(x_{n})=s(x_0)$ and $r(x_k)=s(x_{k+1})$ for the rest indices. The edges $x_k$ here are either of the form $\iota(y_j)$ with $\iota$ as in \eqref{eq:iota}, or belong to the unique directed path in a connected component of $X$ (which is a bi-tree) between $r(\iota(y_j))$ and $s(\iota(y_j+1))$ (resp.\ $r(\iota(y_m))$ and $s(\iota(y_0))$).

Let now $\tilde E$ be a rose-tree bi-graph obtained from $E$ as in Theorem \ref{thm:reduction}. By Lemma \ref{graph_lem_rose_tree_K_1} the elements $u_\mu-\ol u_\mu$, where $u_\mu$ runs through all the loops in ${\tilde E}^{1,0}$ form a basis of the group $\K_1(\Osm_{\tilde E})$. According to the previous paragraph, the preimages of these entities under the corresponding projection map in $\check E$ are of the form $\sum_{j=0}^{m_\mu}(y_{\mu,j}-\ol y_{\mu,j})$ with the properties as above (i.e.\ forming a cycle), and the images of the latter entities look the same in $\check E$. Similarly, the preimages of the established entities $\sum_{j=0}^{m_\mu}(y_{\mu,j} - \ol y_{\mu,j})\in \Zb^{{\check E}^{1}}$ in $\Zb^{E^1}$ are of the form $\sum_{k=0}^{n_\mu}(x_{\mu,k}-\ol x_{\mu,k})$ with similar properties. The union of the edges $x_{\mu,k}$ defines a cycle in the geometric realization of $\E$. These cycles are not homotopic to each other for different $\mu$, since the geometric realizations of $\E$ and $\tilde \E$ are homotopic by the construction and a chain homotopy implemented by the sequence $\E,\hat \E, \check \E,\tilde \E$ transfers these loops into the ones given by $u_\mu$ in the realization $\tilde \E$, which are obviously not homotopic for different $\mu$.

Thus, we have established a map from the set of generators of $\K_1(\Osm_{\tilde \E})$ and distinct loops in the realization of $\E$.
\end{proof}
\end{cor}


\begin{cor}
There is a natural homomorphism $\K_0(\Osm_\E)\cong H_1(\E)\oplus \Zb^{\gamma(\E)}$.
\begin{proof}
Here we only sketch the proof, as the result is obvious at the intuitive level whereas rigorous arguments will take too much place. In order to prove the claim, we again make a closer look onto the proof of Theorem \ref{thm:reduction}.

Indeed, we may trace the spanning bi-trees $\check T_{\mu,2k}$ in $\check E$ back to bi-trees $T_{\mu,2k}$ in the bi-graph $E$ in the unique way. It then follows from the construction that the union of the bi-trees $T=\bigcup_{\mu,k}T_{\mu,k}$ is a spanning bi-tree for the bi-graph $E$.

Consider now the set of bi-loops $\{ u_\nu,\ol u_\nu\}$ in the graph $\tilde E$. As we have seen in Lemma \ref{graph_lem_rose_tree_K_0}, the set $\{u_\nu\}$ may be regarded as a set of independent generators (but still not the whole set of generators) of the group $\K_0(\Osm_{\tilde E})$. For that reason, the elements $u+\Im(\Id_{\check E}-\Phi_{\check E})$ form a set of independent generators in $\K_0(\Osm_{\check E})$ (cf.\ Lemma \ref{graph_lem_getridoftrees_0}). The latter elements, in turn, have unique preimages in $\coker(\Id_{\hat E}-\Phi_{\hat E})$, and it can be checked that their images in $\coker(\Id_{\hat E}-\Phi_{\hat E})$ may be attributed to those (real) edges in de bi-graph $E$ that do not belong to the tree bi-graph $T$ introduced above. In the spirit of Lemma \ref{graph_lem_getridoftrees_0} we abusively denote these elements by $\{u_\nu\}$. To each element $u_\nu$ in $E$ we may attribute a cycle in $\E$ that contains $u_\nu$ and connects $s(u_\nu)$ and $r(u_\nu)$ via the unique path in $T$. These cycles are independent by the same argument as in Corollary \ref{cor:inv:K_1}, and so may be used to define a basis of the group $H_1(\E)$. As now the set $\{u_\nu\}$ constitutes a basis of a subgroup of $\K_0(\Osm_{\tilde E})$, we have established a natural inclusion
\[
H_1(\E)\hookrightarrow \K_0(\Osm_{\tilde \E}) \cong \K_0(\Osm_\E).
\]

As for the rest of the generators, we have observed in Subsection \ref{subsec:Vanency} that there is a natural isomorphism between the valency sets of two graphs, one of which is obtained from another by shrinking of a countable union of mutually non-intersecting finite subgraphs. Hence, we may establish natural isomorphisms of valency sets $V_{\E}\cong V_{\hat \E} \cong V_{\check \E} \cong V_{\tilde \E}$, which yields an isomorphism between the rest of the generators of $\K_0(\Osm_\E)$ and the valency set of $\E$.

Summarizing, we see that, indeed, $\K_0(\E) = \Zb^{\beta(\E)}\oplus \Zb^{\gamma(\E)}$.

\end{proof}
\end{cor}

\vskip1truecm

{\bf Acknowledgements}

We are grateful to the Max-Planck-Institute for Mathematics in Bonn, where most part of this research have been done, for hospitality, support, and excellent research atmosphere. 
This work is funded by the ERC grant 320974.

\end{document}